\PassOptionsToPackage{square,numbers}{natbib}
\documentclass[10pt]{article}

\usepackage[utf8]{inputenc}
\usepackage[T1]{fontenc}
\usepackage[english]{babel}
\usepackage{graphicx}
\usepackage{amsmath,amssymb,amsbsy,amsthm,color,natbib}
\usepackage{multicol}
\usepackage{color} 
\usepackage{natbib}
\usepackage[linktoc=all,colorlinks=true,citecolor=blue,linkcolor=blue]{hyperref}
\usepackage[paperwidth=18cm,paperheight=29.7cm,hmargin=2cm,vmargin=3cm]{geometry}
\usepackage{authblk}
\usepackage{comment}

\title{\textbf{Chaos and Regularity in the Double Pendulum \\ with Lagrangian Descriptors}}

\author[1]{Javier Jiménez López\thanks{javier.jimenezl@edu.uah.es}}
\author[1]{V. J. García-Garrido\thanks{ vjose.garcia@uah.es (Coresponding author).}}
\affil[1]{Departamento de F\'isica y Matem\'aticas, Facultad de Ciencias, Universidad de Alcal\'a, 28805 Alcal\'a de Henares, Madrid, Spain.}

\begin{document}

\maketitle


\begin{abstract}

In this paper we apply the method of Lagrangian descriptors as an indicator to study the chaotic and regular behavior of trajectories in the phase space of the classical double pendulum system. In order to successfully quantify the degree of chaos with this tool, we first derive Hamilton's equations of motion for the problem in non-dimensional form, showing that they can be written compactly using matrix algebra. Once the dynamical equations are obtained, we carry out a parametric study in terms of the system's total energy and the other model parameters (lengths and masses of the pendulums, and gravity), to determine the extent of the chaotic and regular regions in the phase space. Our numerical results show that for a given mass ratio, the maximum chaotic fraction of phase space trajectories is attained when the pendulums have equal lengths. Moreover, we give a characterization of the growth and decay of chaos in the system in terms of the model parameters, and explore the hypothesis that the chaotic fraction follows an exponential law over different energy regimes.

\end{abstract}

\noindent \textbf{keywords:} Hamiltonian systems, Double pendulum, Chaos indicators, Lagrangian descriptors, Chaos and Regularity.


\section{Introduction}

The study of dynamical systems and their chaotic behavior has captured the attention of the scientific community for decades given its prevalence in natural phenomena and its relevance in technological applications \cite{thompson16}. Among these systems, the double pendulum stands out as a classical and fascinating model, not only for its conceptual simplicity but also for the complexity of its dynamical behavior. This system is a perfect example of a nonlinear system, which, despite its apparent mechanical simplicity, exhibits a great variety of dynamical behaviors ranging from periodic motions to chaos.

The purpose of this study is to analyze the fraction of chaotic trajectories in the phase space of the double pendulum, which is a paradigmatic benchmark example of chaos in the nonlinear sciences literature \cite{korsch2007chaos}. In order to carry out this task, we first write the model equations in dimensionless form, and provide a simplified expression for them using matrix algebra. We follow a similar approach as that developed for other mechanical systems such as, for example, the Furuta pendulum \cite{munoz2007bifurcation}. This formulation allows for a more detailed exploration of the degree of chaos present in the system depending on its characteristic parameters. In a recent article concerned with coupled pendulums \cite{szuminski2024new}, the authors highlight the importance of carrying out complete parametric studies for such systems, by means of applying different techniques, as a fundamental step towards obtaining relevant insights about their rich and complex dynamical behavior. In the problem we address in this paper, we use chaos indicators derived from the method of Lagrangian descriptors \cite{Hille22,zimper23} to quantify chaos in the double pendulum. This trajectory-based diagnostic tool has proved to be particularly effective for the analysis of regular and chaotic dynamics, and also for revealing the underlying phase space structures, in many physical systems \cite{Garcia16,naik20,darwish21,raffa23}. This methodology, which is straightforward to implement since it only requires the integration of ensembles of trajectories, has the capability of reducing the computational time required to study the invariant manifolds that divide the phase space into regions where trajectories display qualitatively distinct dynamical behaviors.

The relevance of our study lies not only in its contribution to the understanding of the double pendulum, but also to the analysis of similar Hamiltonian systems such as those arising in the motion of celestial bodies \cite{kaheman2023saddle}, which share fundamental characteristics with the double pendulum. These connections make the findings of our study potentially applicable in a broader context. In the existing literature, several approaches have been undertaken to study chaos in the double pendulum, ranging from experimental works as those presented in \cite{levien1993double,shinbrot1992chaos} to numerical studies such as \cite{Indiati_2016,STACHOWIAK2006417,miyamoto2013long,deleanu2011dynamics}. However, the application of Lagrangian descriptors in this problem offers a new perspective that combines both mathematical rigor with physical intuition. Recently, a preliminary study of the chaotic fraction in the phase space of the double pendulum has been performed in \cite{cabrera2023regular} by means of the Maximum Lyapunov Exponent. However, in this paper the authors only address the case where the masses and lengths of the pendulums are equal to one, and gravity is taken as that of the Earth. Our work goes further and extends the results obtained in \cite{cabrera2023regular}. We do so by first making the double pendulum model dimensionless and use matrix algebra to simplify the equations of motion. Once Hamilton's equations are written in terms of the mass and length ratios, we then conduct a detailed parametric study to determine how the masses and lengths of the pendulums influence the degree of chaos that the system exhibits.

This article is outlined as follows. In Section \ref{Formulation} we develop the non-dimensional Hamiltonian model of the double pendulum and investigate the critical points of its potential energy surface. Section \ref{Methodology} is devoted to the methodology of this work. We explain how one can define different chaos indicators derived from Lagrangian descriptors in order to study the chaotic fraction of the phase space as a function of the energy and the model parameters. In Section \ref{Results} we discuss and interpret the main results obtained from our analysis of the double pendulum system. We finish this work with Section \ref{Conclusions}, where we summarize the conclusions and provide some ideas for future research.

\section{Hamiltonian Formulation of the Mechanical System} \label{Formulation}

We begin this section by deriving the Hamiltonian formulation and Hamilton's equations of motion for the double pendulum. Despite this is a well-known academic exercise, carried out in many courses and textbooks on Analytical Mechanics, see e.g. \cite{goldstein}, we explain here how to obtain a compact and dimensionless expression of the problem, written in terms of matrix algebra. This approach is inspired by the work developed in \cite{munoz2007bifurcation} for the Furuta pendulum, which has similarities with our model system.

\begin{figure}[htbp]
	\centering
	\includegraphics[scale = 0.8]{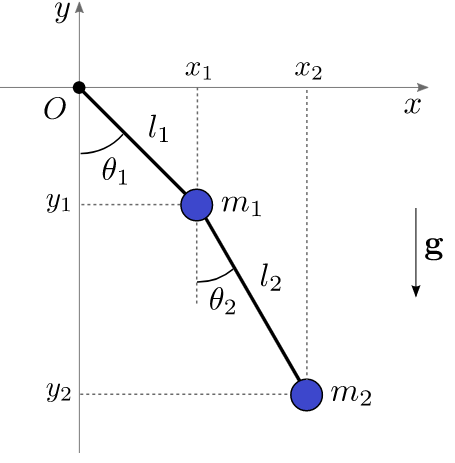}
	\caption{Diagram of the double pendulum system, taken from \cite{AssencioHam}}
	\label{dpend_diag}
\end{figure}

The first step is to define the configuration space vector $\boldsymbol{\theta} = (\theta_1,\theta_2)^T$ that contains the degrees of freedom of the mechanical system. If we refer to the diagram shown in Fig. \ref{dpend_diag}, we have that the position vectors of the masses $m_1$ and $m_2$ for the double pendulum system are given by:
\begin{equation}
\mathbf{r}_1 = l_1\begin{bmatrix}
	\sin \theta_1 \\
	-\cos \theta_1
\end{bmatrix} \quad,\quad \mathbf{r}_2 = \mathbf{r}_1 + \mathbf{r}_{12} = l_1\begin{bmatrix}
\sin \theta_1 \\
-\cos \theta_1
\end{bmatrix} + l_2\begin{bmatrix}
\sin \theta_2 \\
-\cos \theta_2
\end{bmatrix} \;,
\end{equation}
and the velocity vectors:
\begin{equation}
\dot{\mathbf{r}}_1 = l_1 \dot{\theta}_1 \begin{bmatrix}
	\cos \theta_1 \\
	\sin \theta_1
\end{bmatrix} \quad,\quad \dot{\mathbf{r}}_2 = \dot{\mathbf{r}}_1 + \dot{\mathbf{r}}_{12} =  l_1 \dot{\theta}_1 \begin{bmatrix}
\cos \theta_1 \\
\sin \theta_1
\end{bmatrix} + l_2 \dot{\theta}_2 \begin{bmatrix}
\cos \theta_2 \\
\sin \theta_2
\end{bmatrix} \;.
\end{equation}
The kinetic energy of the system is:
\begin{equation}
T = \dfrac{1}{2}m_1 \lVert \dot{\mathbf{r}}_1 \rVert^2 + \dfrac{1}{2}m_2 \lVert \dot{\mathbf{r}}_2 \rVert^2 \;.
\end{equation}
Now, we have that:
\begin{equation}
\lVert \dot{\mathbf{r}}_1 \rVert^2 = \dot{\mathbf{r}}_1 \cdot \dot{\mathbf{r}}_1 = l_1^2 \dot{\theta_1}^2 = \dot{\boldsymbol{\theta}}^T \begin{bmatrix}
	l_1^2 & 0 \\
	0 & 0
\end{bmatrix} \dot{\boldsymbol{\theta}} \;,
\end{equation}
where $\dot{\boldsymbol{\theta}} = (\dot{\theta}_1,\dot{\theta}_2)^T$ is the vector of generalized velocities. The degrees of freedom of the system are $\theta_1$ and $\theta_2$. Moreover:
\begin{equation}
\lVert \dot{\mathbf{r}}_2 \rVert^2 = \dot{\mathbf{r}}_{12} \cdot \dot{\mathbf{r}}_{12} = \lVert \dot{\mathbf{r}}_1 \rVert^2 + \lVert \dot{\mathbf{r}}_{12} \rVert^2 + 2 \dot{\mathbf{r}}_1 \cdot \dot{\mathbf{r}}_{12} \;,
\end{equation}
and using the expressions:
\begin{equation}
\lVert \dot{\mathbf{r}}_{12} \rVert^2 = \dot{\mathbf{r}}_{12} \cdot \dot{\mathbf{r}}_{12} = l_2^2 \dot{\theta_2}^2 = \dot{\boldsymbol{\theta}}^T \! \begin{bmatrix}
	0 & 0 \\
	0 & l_2^2
\end{bmatrix} \dot{\boldsymbol{\theta}} \;,
\end{equation} 
\begin{equation}
2 \dot{\mathbf{r}}_1 \cdot \dot{\mathbf{r}}_{12} = 2 l_1 l_2 \dot{\theta}_1 \dot{\theta}_2 \cos\Delta \theta = \dot{\boldsymbol{\theta}}^T \! \begin{bmatrix}
	0 & l_1 l_2 \cos\Delta \theta \\
	l_1 l_2 \cos\Delta \theta & 0
\end{bmatrix} \dot{\boldsymbol{\theta}} \;,
\end{equation}
where $\Delta \theta = \theta_1 - \theta_2$, we arrive at:
\begin{equation}
T = \dfrac{1}{2} \dot{\boldsymbol{\theta}}^T \! \begin{bmatrix}
	(m_1+m_2)l_1^2 & m_2l_1 l_2  \cos\Delta \theta \\
	m_2l_1 l_2  \cos\Delta \theta & m_2 l_2^2
\end{bmatrix} \dot{\boldsymbol{\theta}} \;.
\end{equation}
This can be rewritten as:
\begin{equation}
T = \dfrac{1}{2} m_2 l_2^2 \,\dot{\boldsymbol{\theta}}^T \! \begin{bmatrix}
	\beta & \alpha \cos \Delta \theta \\
	\alpha \cos \Delta \theta & 1
\end{bmatrix} \dot{\boldsymbol{\theta}} = \dfrac{1}{2} m_2 l_2^2 \dot{\boldsymbol{\theta}}^T B(\cos \Delta \theta) \dot{\boldsymbol{\theta}} \;,
\end{equation}
where we have defined the matrix:
\begin{equation}
B(x) = \begin{bmatrix}
	\beta & \alpha x \\
	\alpha x & 1
\end{bmatrix} \;,
\label{Bmat}
\end{equation}
and introduced the dimensionless quantities:
\begin{equation}
\alpha = \dfrac{l_1}{l_2} > 0 \quad,\quad \beta = \dfrac{m_1+m_2}{m_2} \alpha^2 = \dfrac{\alpha^2}{\mu} > 0 \quad,\quad \mu = \dfrac{m_2}{m_1+m_2} \;.
\end{equation}
If we set $\sigma = \dfrac{m_1}{m_2}$ then:
\[
\mu = \dfrac{1}{1+\sigma} \quad,\quad \beta = {\alpha^2}(1+\sigma) \;.
\]
Notice that by construction, the model parameters satisfy the conditions $\beta > \alpha^2$, $0 < \mu < 1$ and $0 < \sigma < +\infty$. Using this simplified expression for the kinetic energy, the Lagrangian of the system is:
\begin{equation}
\mathcal{L} = T - V = \dfrac{1}{2} m_2 l_2^2 \, \dot{\boldsymbol{\theta}}^T \! B(\cos \Delta \theta) \,\dot{\boldsymbol{\theta}} + (m_1+m_2) \, l_1g\cos\theta_1 + m_2 g l_2 \cos\theta_2 \;.
\end{equation}
If we non-dimensionalize time by setting $t = \sqrt{\dfrac{l_2}{g}} \, \tau$ then:
\begin{equation}
\dot{\boldsymbol{\theta}} = \dfrac{d \boldsymbol{\theta}}{dt} = \dfrac{d \boldsymbol{\theta}}{d\tau} \dfrac{d \tau}{d t} = \sqrt{\dfrac{g}{l_2}} \dfrac{d \boldsymbol{\theta}}{d\tau} = \sqrt{\dfrac{g}{l_2}} \, \boldsymbol{\theta}^\prime \;.
\end{equation}
This yields the non-dimensional kinetic energy:
\begin{equation}
\mathcal{T} = \dfrac{T}{m_2 l_2 g} = \dfrac{1}{2}  \left(\boldsymbol{\theta}^\prime\right)^T \!\! B(\cos \Delta \theta) \, \boldsymbol{\theta}^\prime \;.
\end{equation}
Similarly, the potential energy of the system is:
\begin{equation}
V = -(m_1+m_2)l_1g\cos\theta_1 - m_2 g l_2 \cos\theta_2 = -m_2 l_2 g \left[\dfrac{\beta}{\alpha} \cos \theta_1 + \cos \theta_2\right] \;,
\end{equation}
and if we make it non-dimensional:
\begin{equation}
\mathcal{V} = \dfrac{V}{m_2 l_2 g} = -\dfrac{\beta}{\alpha} \cos \theta_1 - \cos \theta_2 \;.
\label{pot_dpend}
\end{equation}
Hence, the dimensionless lagrangian is:
\begin{equation}
\mathcal{L} = \dfrac{1}{2} \left(\boldsymbol{\theta}^\prime\right)^T \!\! B(\cos \Delta \theta) \, \boldsymbol{\theta}^\prime + \dfrac{\beta}{\alpha} \cos \theta_1 + \cos \theta_2 \;.
\end{equation}
In order to write the Hamiltonian, we need first the conjugate momenta:
\begin{equation}
\mathbf{p} = \nabla_{\boldsymbol{\theta}^\prime} \mathcal{L} = B(\cos \Delta \theta) \, \boldsymbol{\theta}^\prime = (p_1,p_2)^T \;.
\end{equation}
Therefore:
\begin{equation}
\boldsymbol{\theta}^\prime = B^{-1}(\cos \Delta \theta) \, \mathbf{p} \;,
\end{equation}
where:
\begin{equation}
B^{-1}(x) = \dfrac{1}{\beta-\alpha^2x^2}\begin{bmatrix}
	1 & -\alpha x \\
	-\alpha x & \beta
\end{bmatrix} \;.
\end{equation}
The non-dimensional Hamiltonian is:
\begin{equation}
\mathcal{H}(\boldsymbol{\theta},\mathbf{p}) = \mathcal{T}(\boldsymbol{\theta},\mathbf{p}) + \mathcal{V}\left(\boldsymbol{\theta}\right) = \dfrac{1}{2} \mathbf{p}^T B^{-1}(\cos \Delta \theta) \, \mathbf{p} - \dfrac{\beta}{\alpha} \cos \theta_1 - \cos \theta_2 \;,
\end{equation}
and Hamilton's equations are given by the system:
\begin{equation}
\begin{cases}
	\boldsymbol{\theta}^\prime = \nabla_{\mathbf{p}} \mathcal{H} = \nabla_{\mathbf{p}} \mathcal{T} = B^{-1}(\cos \Delta \theta) \, \mathbf{p} \\[.1cm]
	\mathbf{p}^\prime = -\nabla_{\boldsymbol{\theta}} \mathcal{H} = -\nabla_{\boldsymbol{\theta}} \mathcal{T} - \nabla_{\boldsymbol{\theta}} \mathcal{V} 
\end{cases} \;.
\end{equation}
Now:
\begin{equation}
\nabla_{\boldsymbol{\theta}} \mathcal{V} = \begin{bmatrix}
	\dfrac{\beta}{\alpha} \sin\theta_1 \\[.3cm]
	\sin\theta_2
\end{bmatrix} \;,
\end{equation}
and to calculate $\nabla_{\boldsymbol{\theta}} \mathcal{T}$ we need the derivative of the matrix $B^{-1}$. From the matrix identity:
\begin{equation}
B B^{-1} = \mathbb{I} \;,
\end{equation}
we can differentiate this expression to obtain:
\begin{equation}
\dfrac{dB}{dx} B^{-1} + B \dfrac{d}{dx}\left(B^{-1}\right) = \mathbb{O} \;,
\end{equation}
where $\mathbb{O}$ is the zero matrix. Therefore:
\begin{equation}
\dfrac{d}{dx}\left(B^{-1}\right) = -B^{-1} \dfrac{dB}{dx} B^{-1} \;,
\end{equation}
and taking the derivative of Eq. \eqref{Bmat}:
\begin{equation}
\dfrac{dB}{dx} = \begin{bmatrix}
	0 & \alpha \\
	\alpha & 0
\end{bmatrix} \;.
\end{equation}
This procedure yields the result:
\begin{equation}
C(x) = -B^{-1} \dfrac{dB}{dx} B^{-1} = -\dfrac{\alpha}{\left(\beta-\alpha^2x^2\right)^2} \begin{bmatrix}
	-2\alpha x & \beta+\alpha^2x^2 \\
	\beta+\alpha^2x^2 & -2\alpha\beta x
\end{bmatrix} \;.
\end{equation}
which is a symmetric matrix. We arrive at:
\begin{equation}
\nabla_{\boldsymbol{\theta}} \mathcal{T} = \dfrac{1}{2} \mathbf{p}^T C(\cos \Delta \theta) \, \mathbf{p} \, \nabla_{\boldsymbol{\theta}}\left(\cos \Delta \theta \right) \quad,\quad \nabla_{\boldsymbol{\theta}}\left(\cos \Delta \theta \right) = -\sin \Delta\theta \; \nabla_{\boldsymbol{\theta}}\left( \Delta \theta \right) \;,
\end{equation}
where $\nabla_{\boldsymbol{\theta}}\left( \Delta \theta \right) = (1,-1)^T$. The resulting equations of motion are:
\begin{equation}
\begin{cases}
    \boldsymbol{\theta}^\prime = B^{-1}(\cos \Delta \theta) \, \mathbf{p} \\ 
    \mathbf{p}^\prime = \dfrac{\sin \Delta \theta}{2} \mathbf{p}^T C(\cos \Delta \theta) \, \mathbf{p} \begin{bmatrix}
		1 \\
		-1
	\end{bmatrix} - \begin{bmatrix}
	\dfrac{\beta}{\alpha} \sin\theta_1 \\[.1cm]
	\sin\theta_2
	\end{bmatrix}
\end{cases} \;.
\label{ham_eqs}
\end{equation}
If we introduce at this point the state vector $\boldsymbol{\sigma} = (\boldsymbol{\theta},\mathbf{p})^T$, Hamilton's equations can be written in compact form as:
\begin{equation}
    \boldsymbol{\sigma}^\prime = \mathbf{F}(\boldsymbol{\sigma}) \;,
\end{equation}
where $\mathbf{F} = (\mathbf{F}_1,\mathbf{F}_2)^T$ is the Hamiltonian vector field whose components are defined by:
\begin{equation}
    \mathbf{F}_1(\boldsymbol{\sigma}) = B^{-1}(\cos \Delta \theta) \, \mathbf{p} \quad,\quad \mathbf{F}_2(\boldsymbol{\sigma}) = \dfrac{\sin \Delta \theta}{2} \mathbf{p}^T C(\cos \Delta \theta) \, \mathbf{p} \begin{bmatrix}
		1 \\
		-1
	\end{bmatrix} - \begin{bmatrix}
	\dfrac{\beta}{\alpha} \sin\theta_1 \\[.2cm]
	\sin\theta_2
	\end{bmatrix} \;.
\end{equation}
The Jacobian matrix determining the variational equations of the linearized dynamics about the phase space point $\boldsymbol{\sigma}_0 = (\boldsymbol{\theta}_0,\mathbf{p}_0)^T$ has the block form:
\begin{equation}
    D_{\boldsymbol{\sigma}}\mathbf{F}(\boldsymbol{\sigma}_0) = \begin{bmatrix}
         D_{\boldsymbol{\theta}}\mathbf{F}_1(\boldsymbol{\sigma}_0) & D_{\mathbf{p}}\mathbf{F}_1(\boldsymbol{\sigma}_0) \\[.2cm]
         D_{\boldsymbol{\theta}}\mathbf{F}_2(\boldsymbol{\sigma}_0) & D_{\mathbf{p}}\mathbf{F}_2(\boldsymbol{\sigma}_0)
    \end{bmatrix} \;,
\end{equation}
where the subscript indicates the variables with respect to which we differentiate to obtain each of the matrices involved in the jacobian. After some algebra, the resulting matrix is:
\begin{equation}
    D_{\boldsymbol{\sigma}}\mathbf{F} = \begin{bmatrix}
        -\sin \Delta\theta \; C\left(\cos \Delta \theta\right) \left(\mathbf{p} \otimes \nabla_{\boldsymbol{\theta}}\left(\Delta \theta \right)\right) & B^{-1}(\cos \Delta \theta) \\[.2cm]
        G(\boldsymbol{\theta},\mathbf{p})-H_{\mathcal{V}} & \sin \Delta \theta \left(\nabla_{\boldsymbol{\theta}}\left(\Delta \theta \right) \otimes \mathbf{p}\right) C\left(\cos \Delta \theta\right)
    \end{bmatrix} \;,
\end{equation}
such that $G(\boldsymbol{\theta},\mathbf{0}) = \mathbb{O}$
and the hessian matrix of the potential energy is:
\begin{equation}
    H_{\mathcal{V}}(\boldsymbol{\theta}) = \begin{bmatrix}
        \dfrac{\beta}{\alpha} \cos\theta_1 & 0 \\[.1cm]
        0 & \cos\theta_2 \\[.1cm]
    \end{bmatrix} \;.
\end{equation}

\begin{figure}[htbp]
	\centering
	\includegraphics[scale = 0.29]{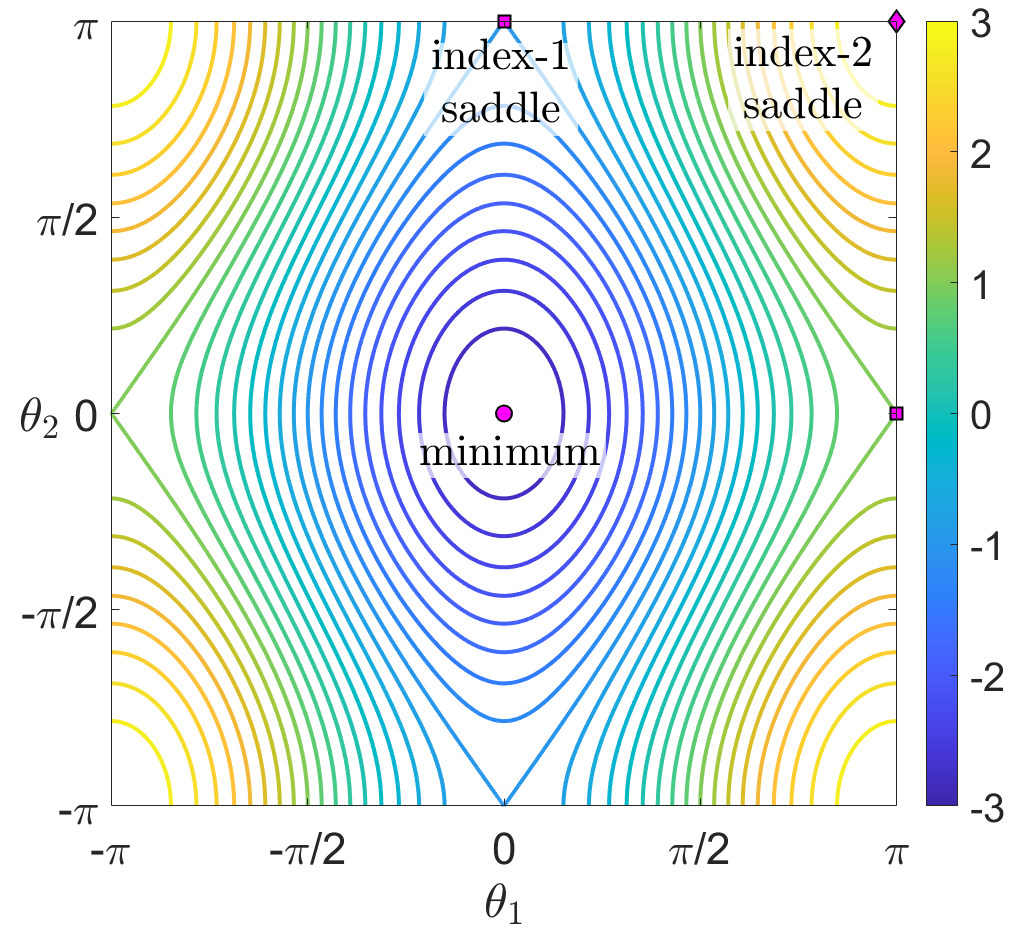}
    \caption{Contours of the potential energy function given in \eqref{pot_dpend} for the double pendulum system with $\alpha = 1$ and $\beta = 2$. We have marked the equilibrium points of Hamilton's equations using magenta squares (saddles), circles (local minimum) and diamonds (local maximum).}
    \label{pes_dpend}
\end{figure}

The dynamical system in Eq. \eqref{ham_eqs} has 4 equilibria, denoted by $\boldsymbol{\sigma}_{\text{eq}} = (\theta_{1,\text{eq}},\theta_{2,\text{eq}},0,0)$, located in configuration space ($\mathbf{p} = \mathbf{0}$) at the points:
\begin{equation}
\theta_{1,\text{eq}} = 0, \pi \quad,\quad \theta_{2,\text{eq}} = 0, \pi  \;.
\end{equation}
and we will denote them as follows:
\begin{equation}
    \boldsymbol{\sigma}^1_{\text{eq}} = (0,0,0,0) \;\;,\;\; \boldsymbol{\sigma}^2_{\text{eq}} = (0,\pi,0,0) \;\;,\;\; \boldsymbol{\sigma}^3_{\text{eq}} = (\pi,0,0,0)  \;\;,\;\; \boldsymbol{\sigma}^4_{\text{eq}} = (\pi,\pi,0,0) \;,
\end{equation}
These points are critical points of the potential energy in Eq. \eqref{pot_dpend}, and in order to get a better understanding of the geometrical features of the potential energy surface, we have represented the function in Fig. \ref{pes_dpend}, together with the equilibria of Hamilton's equations. The energies of the equilibrium points are given by:
\begin{equation}
\begin{cases}
\mathcal{H}_1 = \mathcal{V}_1(0,0) = -\alpha(1+\sigma)-1 \\[.3cm]
\mathcal{H}_2 = \mathcal{V}_2(0,\pi) = -\alpha(1+\sigma)+1 \\[.3cm]
\mathcal{H}_3 = \mathcal{V}_3(\pi,0) = \alpha(1+\sigma)-1 \\[.3cm]
\mathcal{H}_4 = \mathcal{V}_4(\pi,\pi) = \alpha(1+\sigma)+1 
\end{cases} \;.
\end{equation}
Notice that the energy levels of the equilibrium points of the system satisfy the following relations, depending on the values of the model parameters. If $\alpha(1+\sigma) > 1$ then:
\begin{equation}
\mathcal{H}_1 < \mathcal{H}_2 < \mathcal{H}_3 < \mathcal{H}_4 \;,
\end{equation}
whereas in the case $\alpha(1+\sigma) < 1$ we have:
\begin{equation}
\mathcal{H}_1 < \mathcal{H}_3 < \mathcal{H}_2 < \mathcal{H}_4 \;.
\end{equation}
In particular, there is a critical value of the parameters, $\alpha(1+\sigma) = 1$ for which $\mathcal{H}_2 = \mathcal{H}_3$. This means that, for $\beta = \alpha$, the saddles located at the configuration space points $(0,\pi)$ and $(\pi,0)$ have the same energy. 

Next, we will carry out a stability analysis of the equilibria, with the purpose of determining the local behavior of trajectories about these points. First we define the constant $\gamma = \beta-\alpha^2$, and calculate the Jacobian matrix evaluated at each equilibrium point of the system, which yields:
\begin{equation}
    D_{\boldsymbol{\sigma}}\mathbf{F}(\boldsymbol{\sigma}^1_{\text{eq}}) = \dfrac{1}{\gamma}\begin{bmatrix}
        0 & 0 & 1 & -\alpha \\[.1cm]
        0 & 0 & -\alpha & \beta \\[.1cm]
        -\dfrac{\beta\gamma}{\alpha} & 0 & 0 & 0 \\[.25cm]
        0 & -\gamma & 0 & 0
    \end{bmatrix} \quad,\quad D_{\boldsymbol{\sigma}}\mathbf{F}(\boldsymbol{\sigma}^2_{\text{eq}}) = \dfrac{1}{\gamma}\begin{bmatrix}
        0 & 0 & 1 & \alpha \\[.1cm]
        0 & 0 & \alpha & \beta \\[.1cm]
        -\dfrac{\beta\gamma}{\alpha} & 0 & 0 & 0 \\[.25cm]
        0 & \gamma & 0 & 0
    \end{bmatrix} \;,
\end{equation}
\begin{equation}
    D_{\boldsymbol{\sigma}}\mathbf{F}(\boldsymbol{\sigma}^3_{\text{eq}}) = \dfrac{1}{\gamma} \begin{bmatrix}
        0 & 0 & 1 & \alpha \\[.1cm]
        0 & 0 & \alpha & \beta \\[.1cm]
        \dfrac{\beta\gamma}{\alpha} & 0 & 0 & 0 \\[.25cm]
        0 & -\gamma & 0 & 0
    \end{bmatrix} \quad,\quad D_{\boldsymbol{\sigma}}\mathbf{F}(\boldsymbol{\sigma}^4_{\text{eq}}) = \dfrac{1}{\gamma} \begin{bmatrix}
        0 & 0 & 1 & -\alpha \\[.1cm]
        0 & 0 & -\alpha & \beta \\[.1cm]
        \dfrac{\beta\gamma}{\alpha} & 0 & 0 & 0 \\[.25cm]
        0 & \gamma & 0 & 0 
    \end{bmatrix} \;.
\end{equation}
Analyzing the eigenvalues of these matrices, it can be shown that $\boldsymbol{\sigma}^1_{\text{eq}}$ is stable, with center$\times$center stability, $\boldsymbol{\sigma}^2_{\text{eq}}$ and $\boldsymbol{\sigma}^3_{\text{eq}}$ are index-1 saddles, that is, they have saddle$\times$center stability. Finally, $\boldsymbol{\sigma}^4_{\text{eq}}$ is an index-2 saddle, with saddle$\times$saddle stability.


\section{Methodology} \label{Methodology}

In this section we describe the methodology that we have followed to quantify the chaotic and regular fraction of the phase space in the double pendulum system. There exists many chaos indicators in the literature to carry out this task \cite{skokos2016chaos}, however they typically require the solution of the variational equations of motion together with Hamilton's equations in order to classify a trajectory as chaotic or regular. Moreover, the integration time for applying these techniques is typically of the order $10^6$, which would require symplectic integrators to preserve the energy, or the use of high order Runge-Kutta numerical schemes. Here we explore the use of the method of Lagrangian descriptors (LDs) to carry out the classification of ensembles of trajectories. This new approach has been validated recently on continuous and discrete dynamical systems by developing several chaos indicators based on the values of the LD function \cite{Hille22,zimper23}. These diagnostics have been benchmarked against the Small Alignment Index method and have shown a success rate for the classification of trajectories of above $90\%$. Moreover, the LD chaos indicators only require an integration time of about three orders of magnitude less than the classical methods. Furthermore, the implementation of LDs is straightforward, since one only needs to add an extra differential equation to the system of Eqs. \ref{ham_eqs}. This makes it a suitable tool for all practitioners of Dynamical Systems Theory and offers an important advantage with respect to computational time.

The method of Lagrangian descriptors is a trajectory-based diagnostic technique that was originally developed in the field of Geophysics to analyze lagrangian transport and mixing processes in the ocean and the atmosphere \cite{madrid2009,mendoza10}. A lagrangian descriptor is a scalar function constructed in the following way. Given a dynamical system in the form:
\begin{equation}
    \dfrac{d\mathbf{x}}{dt} = \mathbf{f}(\mathbf{x},t) \;.
    \label{general_ds}
\end{equation}
we define a non-negative function $\mathcal{F}(\mathbf{x}(t;\mathbf{x}_0),t)$ that depends on the initial condition $\mathbf{x}_0$ at time $t=t_0$. To determine the LD scalar field, denoted by $\mathcal{L}$, we set an integration time $\tau > 0$ and calculate:
\begin{equation}
    \mathcal{L}(\mathbf{x}_0,t_0,\tau) = \mathcal{L}_f(\mathbf{x}_0,t_0,\tau) + \mathcal{L}_b(\mathbf{x}_0,t_0,\tau) \;,
\end{equation}
where the forward ($\mathcal{L}_f$) and backward ($\mathcal{L}_b$) components of the LD function are given by:
\begin{equation}
\mathcal{L}_f(\mathbf{x}_0,t_0,\tau) = \int_{t_0}^{t_0+\tau} \mathcal{F}(\mathbf{x}(t;\mathbf{x}_0),t) \, dt \quad,\quad
\mathcal{L}_b(\mathbf{x}_0,t_0,\tau) = \int_{t_0-\tau}^{t_0} \mathcal{F}(\mathbf{x}(t;\mathbf{x}_0),t) \, dt \;. 
\end{equation}
Hence, the calculation of LDs involve accumulating the value of the function $\mathcal{F}$ along the trajectory starting at $\mathbf{x}_0$ as it evolves forward and backward in time. In the literature it has been shown rigorously that the scalar field generated by the method has the capability of identifying the invariant sets (equilibria, stable and unstable manifolds, tori, periodic orbits, etc,) that characterize the dynamical behavior of trajectories in the phase space of the system \cite{mancho2013,lopesino2017}. 

\begin{figure}[htbp]
    \begin{center}
    A)\includegraphics[scale = 0.3]{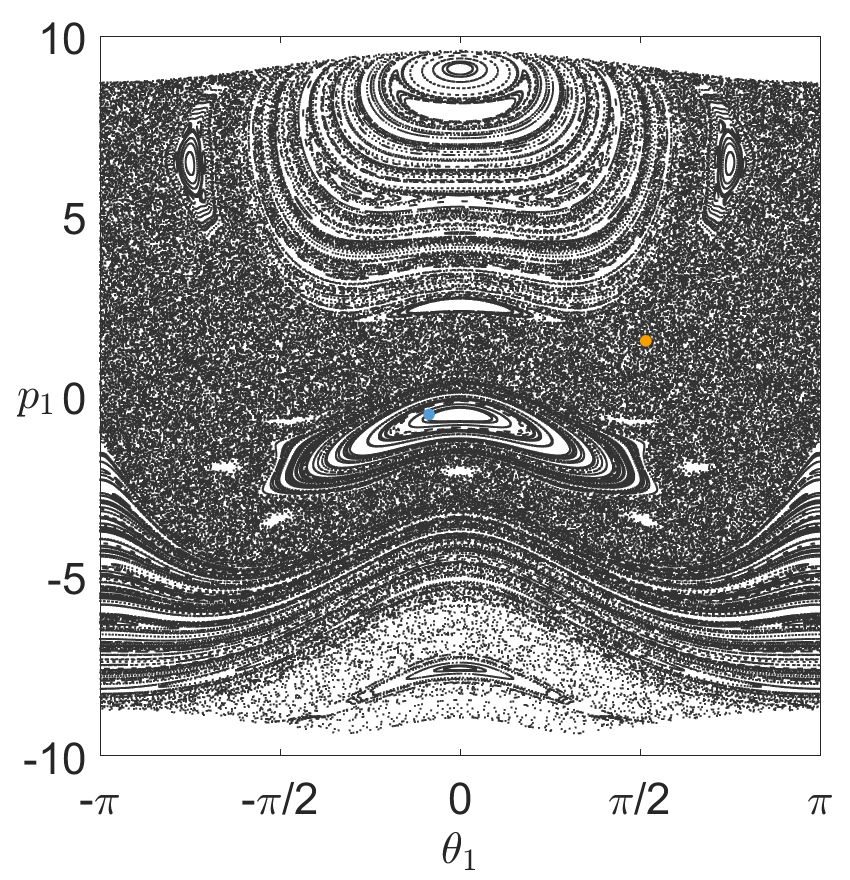}
    B)\includegraphics[scale = 0.61]{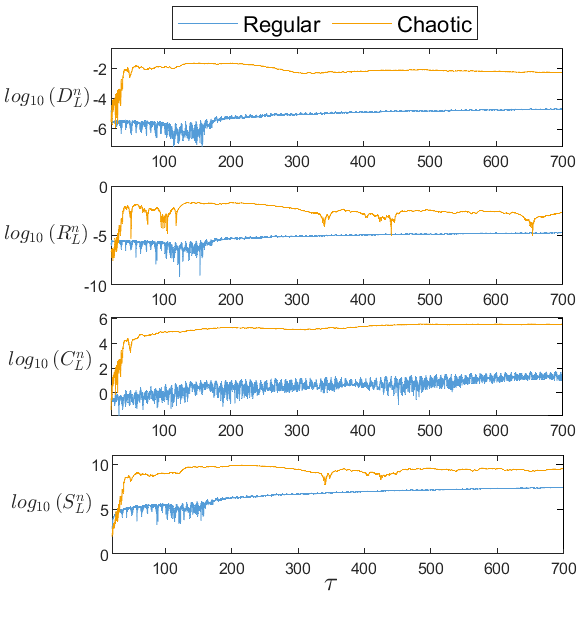}
    \end{center}
    \caption{A) Poincaré section in the phase space slice given by Eq. \eqref{psec} for the system with $\alpha = 1$, $\sigma = 1$ and $\mathcal{H}_0 = 20$. We have marked with an orange point a chaotic initial condition, and with blue a regular one; B) Chaos indicators defined in Eq. \eqref{chaos_inds} calculated for the initial conditions in panel A).}
    \label{Psec_chaos_reg}
\end{figure}

In this work we will use the LD whose definition resembles the $p$-norm of Functional Analysis \cite{lopesino2017}. Thus, we select as the non-negative function the expression:
\begin{equation}
    \mathcal{F}(\mathbf{x},t) = \sum_{i=1}^n |f_i(\mathbf{x},t)|^p \;\;,\;\; 0 < p \leq 1 \;,
\end{equation}
where $f_i$ is the $i$-th component of the vector field that determines the dynamical system in Eq. \ref{general_ds}. In particular we will set $p = 1/2$, which is commonly used because it yields a scalar field where the phase space structures are nicely detected. Notice that the classification of a trajectory as chaotic or regular is independent of the fact that we integrate it forward or backward in time, since both directions must yield the same answer. Hence, for our analysis we have decided to integrate only forward, focusing on the output provided by the component $\mathcal{L}_f$. For the integration time, we have chosen a value of $\tau = 700$ which is close to the one used in previous works \cite{Hille22,zimper23} for the computation of the LD chaos indicators. This value has been adjusted by looking at several trajectories that are known a priori to be chaotic and regular for the double pendulum, see Fig. \ref{Psec_chaos_reg}. The numerical method that we have used to solve Hamilton's equations is an adaptive Runge-Kutta scheme of order 8 (in particular, the Dormand-Prince 8(9) algorithm) with an error tolerance of $10^{-8}$, which has achieved a conservation of the energy with at least 4 significant digits for all the cases considered in this work. 

In order to study the chaotic and regular fraction of the phase space in the double pendulum system, we will carry out our simulations on the Poincaré section defined as follows:
\begin{equation}
    \Sigma\left(\mathcal{H}_0\right) = \left\lbrace (\theta_1,\theta_2,p_1,p_2) \in \mathbb{R}^4 \; \Big| \; \mathcal{H} = \mathcal{H}_0 \;,\; \theta_2 = 0 \;,\; \beta p_2 - \alpha p_1 \cos(\Delta \theta) > 0  \right\rbrace \;.
    \label{psec}
\end{equation}
Note that, given the symmetry of the system with respect to the angle $\theta_1$, we will only consider positive values for this configuration variable, that is, $\theta_1 \in [0,\pi]$. To use the chaos indicators derived from Lagrangian descriptors, we need to obtain the neighbors of the initial condition we want to analyze. These neighbors are given by:
\begin{equation}
\mathbf{y}_i^{\pm} = \mathbf{x}_0 \pm \sigma_i \, \mathbf{e}_i \;\;,\;\; i=1,\ldots,n \;,
\end{equation}
where $\mathbf{e}_i$ is the $i$-th canonical basis vector in $\mathbb{R}^n$, and $\sigma_i$ represents the distance between the central point $\mathbf{x}_0$ and its neighbors on the grid. The value of $n$ represents the dimension of the space where the initial conditions are selected. In particular, since our phase space slice is a plane, then $n = 2$, and hence each initial condition on the grid has 4 neighbors. For our simulations we have chosen a value of $\sigma_i = 10^{-4}$. Using these points, one can construct different chaos indicators based on LDs. In this paper we use the following four diagnostics:
\begin{equation}
\begin{split}
    D_L^n(\mathbf{x}_0) &= \dfrac{1}{2n \mathcal{L}_{f}(\mathbf{x}_0)} \sum_{i=1}^{n} \big|\mathcal{L}_{f}(\mathbf{x}_0)-\mathcal{L}_{f}\left(\mathbf{y}_i^{+}\right)\big|+\big|\mathcal{L}_{f}(\mathbf{x}_0)-\mathcal{L}_{f}\left(\mathbf{y}_i^{-}\right)\big| \;, \\  
    R_L^n(\mathbf{x}_0) &= \Bigg|1-\dfrac{1}{2n \mathcal{L}_{f}(\mathbf{x}_0)} \sum_{i=1}^{n} \mathcal{L}_{f}\left(\mathbf{y}_i^{+}\right)+\mathcal{L}_{f}\left(\mathbf{y}_i^{-}\right) \Bigg| \;, \\[.1cm]
    C_L^n(\mathbf{x}_0) &= \dfrac{1}{2n} \sum_{i=1}^{n} \dfrac{\big|\mathcal{L}_{f}\left(\mathbf{y}_i^{+}\right)-\mathcal{L}_{f}\left(\mathbf{y}_i^{-}\right)\big|}{\sigma_i} \;, \\[.1cm]
    S_L^n(\mathbf{x}_0) &= \dfrac{1}{n} \sum_{i=1}^{n} \dfrac{\big|\mathcal{L}_{f}\left(\mathbf{y}_i^{+}\right)-2\mathcal{L}_{f}(\mathbf{x}_0)+\mathcal{L}_{f}\left(\mathbf{y}_i^{-}\right)\big|}{\sigma_i^2} 
\end{split}
\label{chaos_inds}
\end{equation}
where $\mathcal{L}_{f}(\cdot)$ is the LD value calculated for an integration time $\tau$. Notice that the chaotic or regular nature of a trajectory is equivalently characterized if the we integrate forward or backward, so only one of the components of the LD function is required to determine the chaos indicators. If we set a tolerance $\varepsilon>0$ (that depends on the dynamical system under study and its parameters, and also on the chaos indicator chosen), we say that the initial condition $\mathbf{x}_0$ gives rise to a chaotic trajectory when the logarithm of the value of the indicator at that point is above $\varepsilon$. Otherwise, we classify the trajectory as regular.

To determine the threshold corresponding to each energy and $\alpha$ and $\sigma$ parameter values, we have developed an algorithm based on the histograms generated using the chaos indicators in Eq. \ref{chaos_inds}. Upon examining the shape of these histograms, we have identified three potential scenarios: the first is where the entire phase space is regular, leading to only one peak in the histogram; the second scenario occurs when the phase space contains both regular and chaotic regions, resulting in two peaks in the histogram, with the left one corresponding to the regular fraction and the right one to the chaotic fraction; and the third scenario is when the entire phase space is chaotic, where we will again observe a single peak in the histogram. In the first and third scenarios, determining the threshold is relatively straightforward, but in the second scenario, we have decided to define the threshold as the value of the chaos indicator corresponding to the lowest column between the two peaks, which is the minimum value between the two maximum values. 

Even though the four chaos indicators we have considered have a very similar performance,
we have noticed that the $R^n_{L}$ indicator is more sensitive to the initial conditions given to our algorithm and disagrees with the other three indicators in some of the cases we have studied, while the other three indicators converge to the same classification in the vast majority of cases, being $D^n_{L}$ and $S^n_{L}$ the two best indicators for the cases we have analyzed.

\begin{figure}[htbp]
    \centering
    A)\includegraphics[scale = 0.175]{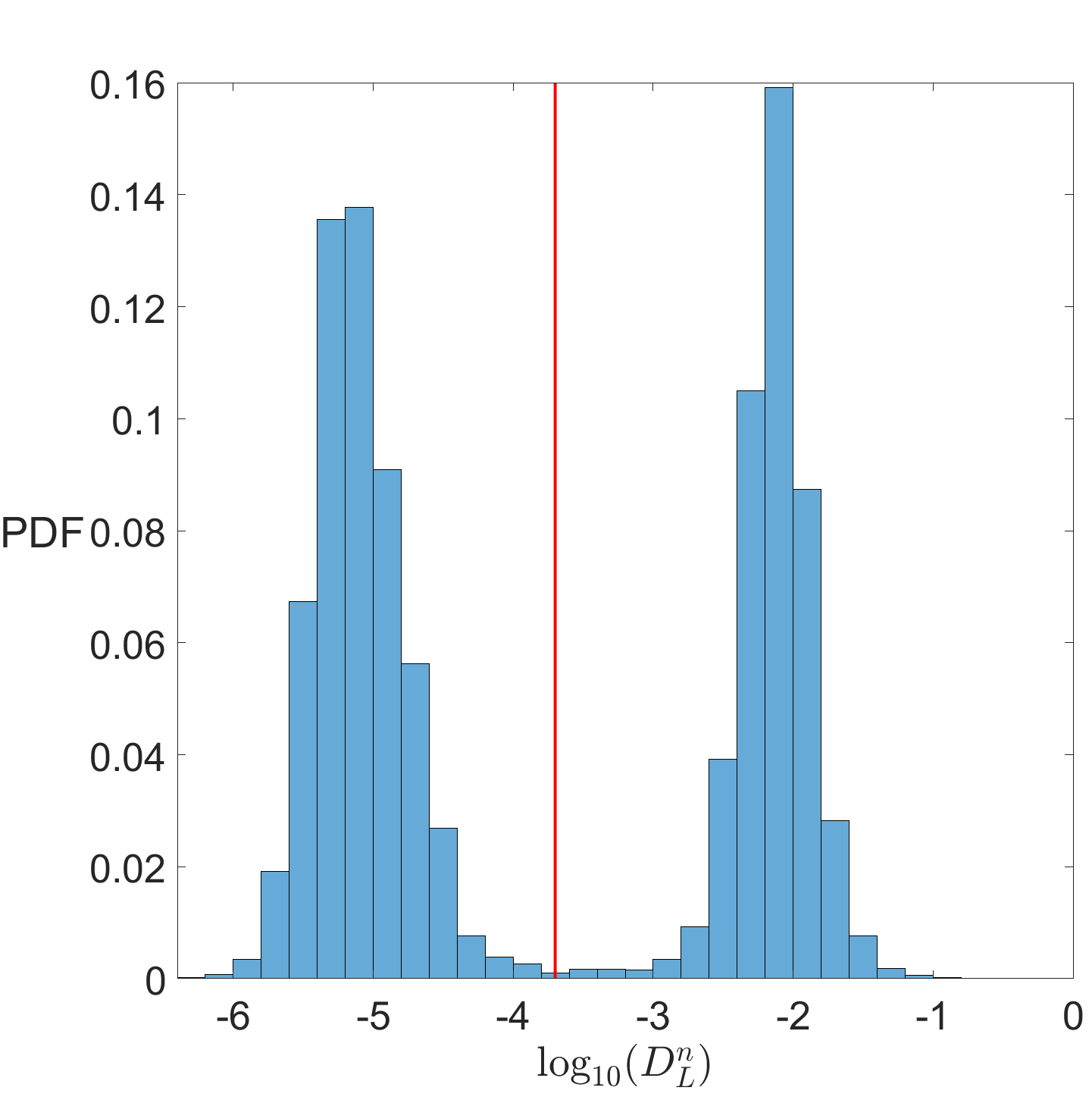}
    B)\includegraphics[scale = 0.435]{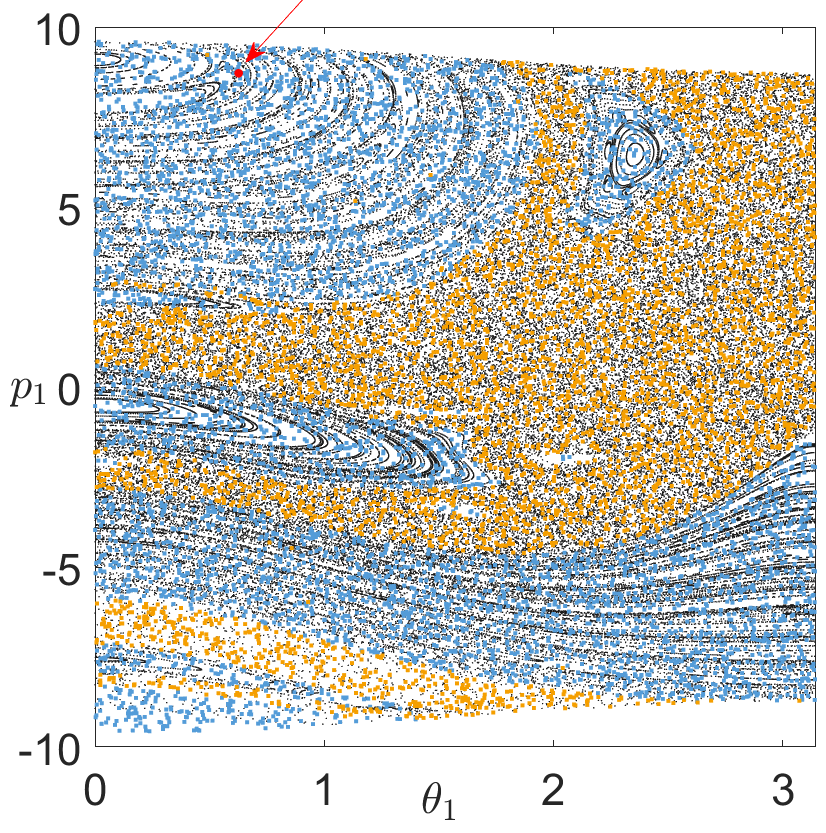} \\
    C)\includegraphics[scale = 0.175]{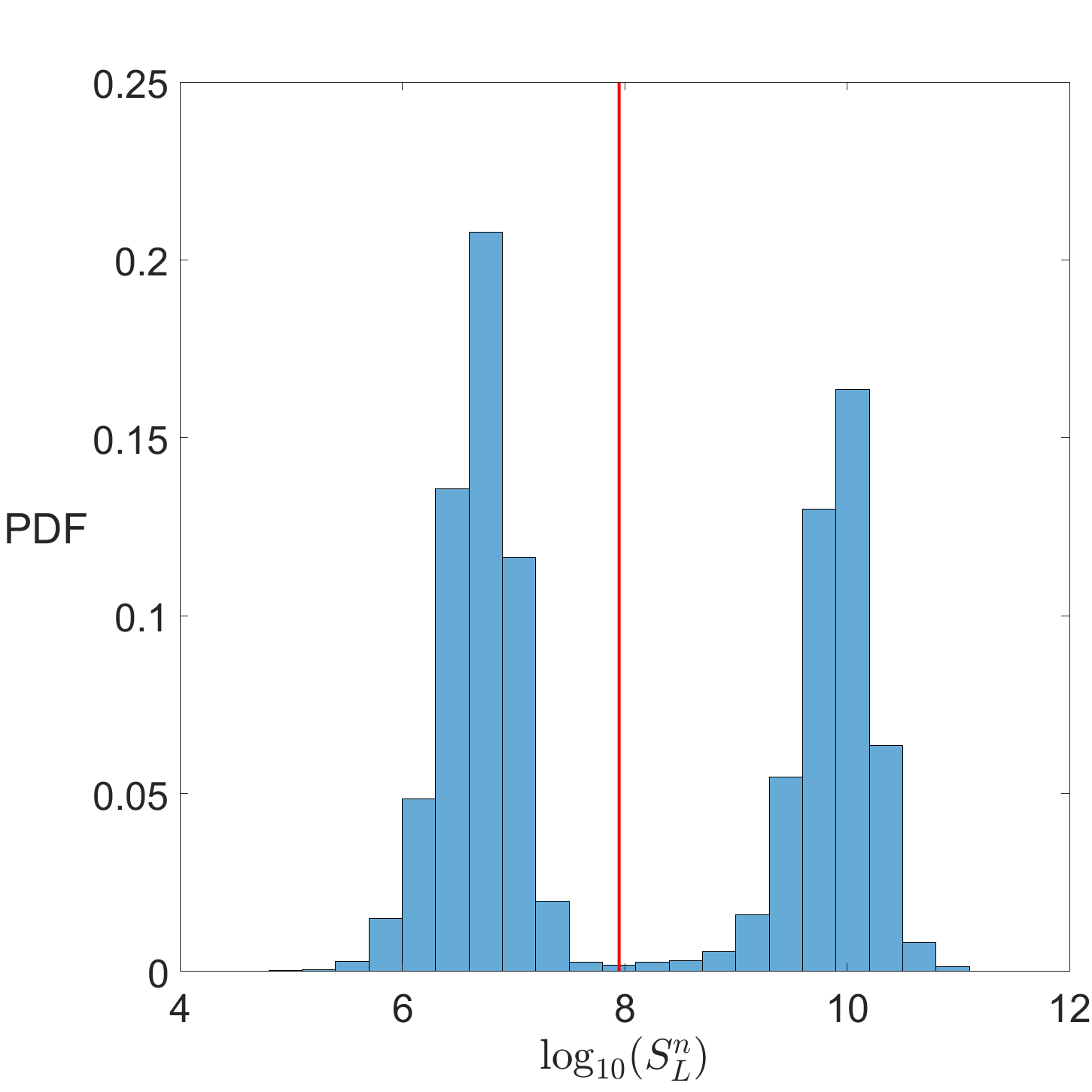} 
    D)\includegraphics[scale = 0.435]{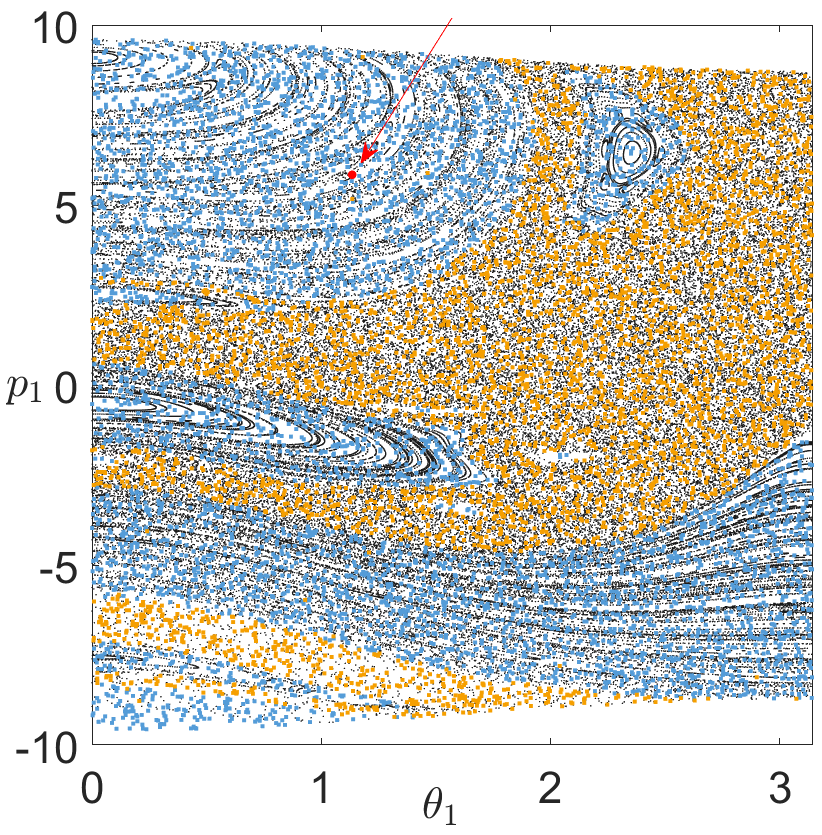}
    \caption{Classification of a random ensemble of $10^4$ initial conditions carried out by means of the chaos indicators based on Lagrangian descriptors. The orange dots correspond to chaotic motion, whereas blue dots indicate regular behavior. For this simulation we have selected an energy of $\mathcal{H}_0 = 20$, and the model parameters are set to $\alpha = 1$ and $\sigma = 1$. A) Histogram showing the distribution of the $D^n_{L}$ indicator. We have marked with a vertical red line the threshold that separates chaotic from regular behavior; B) Poincaré section overlaid with the classified initial conditions using the $D^n_{L}$ indicator. The red dot shows an initial condition that has not been correctly classified by the method.; C) and D) Same analysis, but for the $S^n_{L}$ indicator.}
    \label{S_and_D_panel}
\end{figure}

With this algorithm, we can automate the process of classifying initial conditions in such a way that we can classify any initial condition by comparing the value of the chaos indicator of the condition itself with the threshold we have determined. In Fig. \ref{S_and_D_panel}, the histograms for the $D^n_{L}$ and $S^n_{L}$ indicators are displayed with a vertical red line corresponding to the position of the threshold. From this and following the described process, we have obtained the initial conditions classified by both indicators and we can compare how they have classified the trajectories with the Poincaré section given in Fig. \ref{Psec_chaos_reg} in such a way that it is observed how the chaos indicators have reconstructed the regular and chaotic regions of the phase space as we expected, and very accurately. We want to highlight that, as mentioned at the beginning of this section, the indicators based on Lagrangian descriptors do not classify all the initial conditions correctly. To make this clear, we have indicated in the corresponding figures with dots and red arrows misclassified initial conditions which are regular initial conditions that had been classified as chaotic. Additionally, we have observed that the misclassification is more frequent in the regular region than in the chaotic region of the phase space.

\section{Results} \label{Results}

This section describes all the results we have obtained from the analysis of the chaotic fraction of phase space trajectories in the double pendulum system using the chaos indicators based on LD. In our numerical experiments we have explored the parameter space $(\alpha,\sigma,\mathcal{H}_0)$, where $\alpha = l_1/l_2$ is the ratio of lengths, $\sigma = m_1/m_2$ the ratio of masses and $\mathcal{H}_0$ represents the total energy of the system. For the values of the parameters we have selected $\alpha_i = 2^{i}$ and $\sigma_j=2^{j}$, where $i,j \in \left\lbrace -4,-3,-2,-1,0,1,2,3,4 \right\rbrace$, and we have considered 170 different energy levels for each case $(\alpha_i,\sigma_j)$ of the model parameters. The energies are chosen as follows: 40 energy levels are distributed uniformly from the energy of the local minimum of the potential energy surface ($\mathcal{H}_1$) up to the energy of the local maximum of the potential ($\mathcal{H}_4$), and the other 130 energy values are taken above $\mathcal{H}_4$ with a unit step. For each simulation, we have generated an ensemble of $10^4$ initial conditions, sampled randomly on the surface of section given by Eq. \eqref{psec} using a standard Monte Carlo approach.

\begin{figure}[htbp]
    \centering
    A) \includegraphics[scale=0.21]{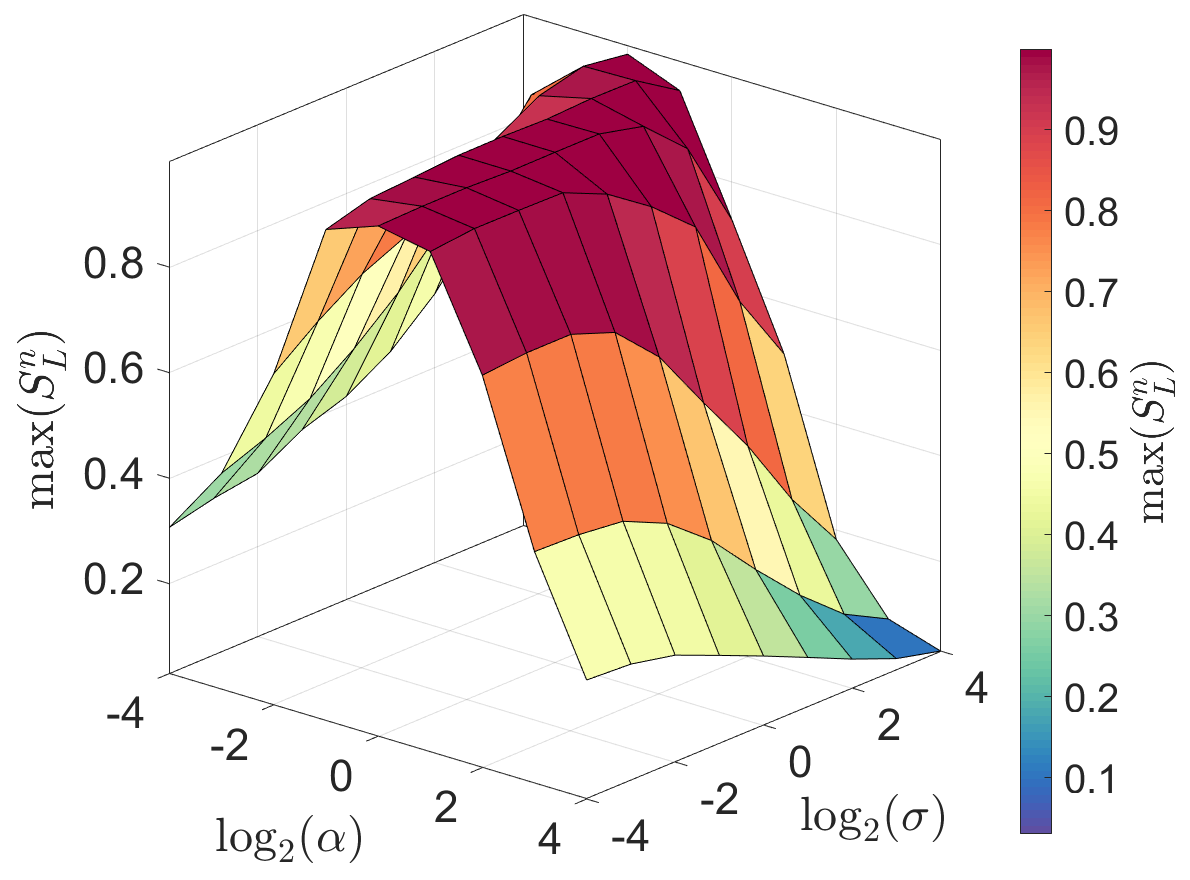} 
    B) \includegraphics[scale=0.215]{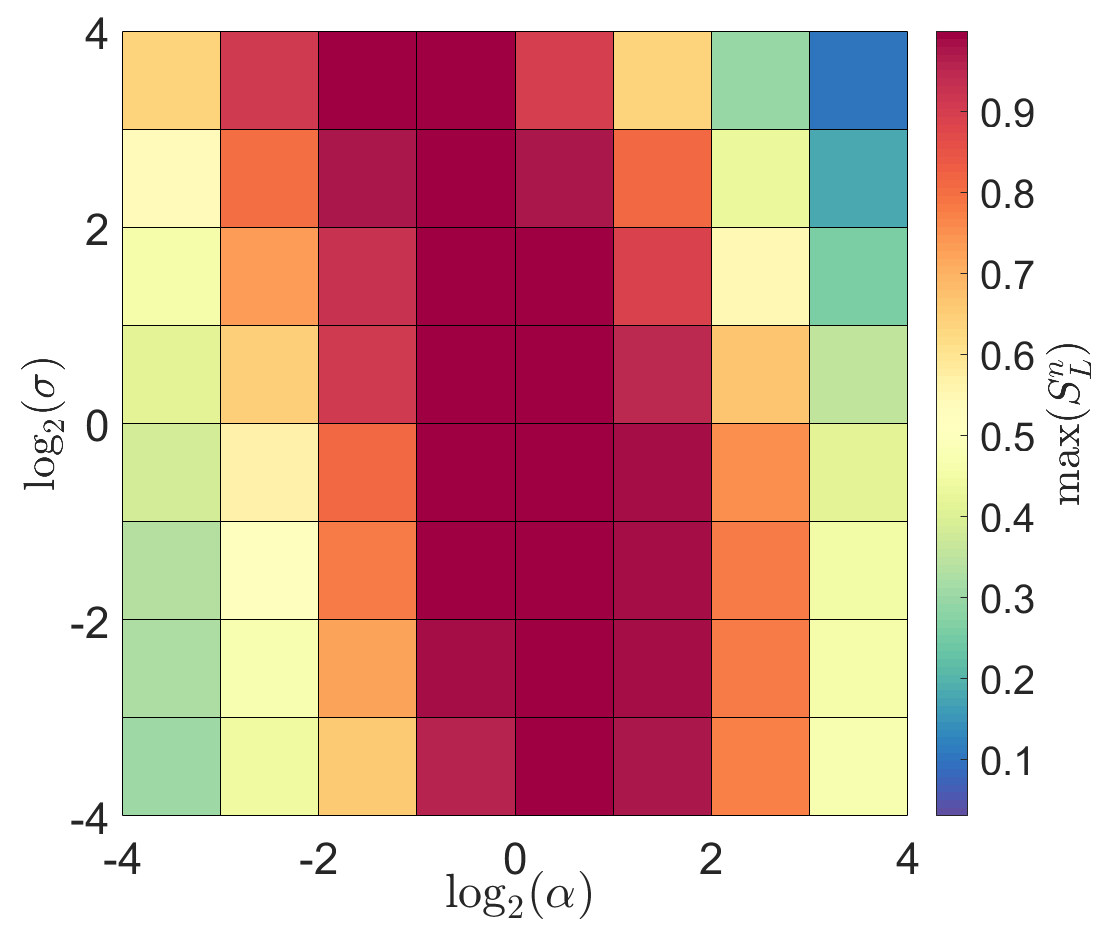}
    \caption{Maximum chaotic fraction of phase space as a function of the model parameters $\alpha = l_1/l_2$ and $\sigma = m_1/m_2$ calculated with the indicator $S^n_L$. The range of energy levels used for the simulations is $[\mathcal{H}_1,\mathcal{H}_4+130]$.}
    \label{max_chaos}
\end{figure}

One of the main results that can be obtained from our simulations is that, for a given value of $\sigma$ (the mass fraction of $m_1$ with respect to $m_2$), the maximum chaotic fraction of the phase space is attained when the lengths of the pendulums are approximately equal, and decreases monotonically as $\alpha$ (the length ratio of $l_1$ with respect to $l_2$) moves away from unity. On the other hand, for a given $\alpha > 1$, chaos increases in the system as $\sigma$ gets smaller, whereas for a fixed $\alpha < 1$, chaos increases with $\sigma$. In particular, this means that if the first pendulum is longer than the second, chaos is enhanced by making the second mass heavier than the first one. On the contrary, in the case where the second pendulum is longer, chaos grows as the mass of the first pendulum gets larger than that of the second pendulum. All this information about the behavior of the system is summarized in Fig. \ref{max_chaos}.

\begin{figure}[htbp]
    \centering
    \includegraphics[scale=0.31]{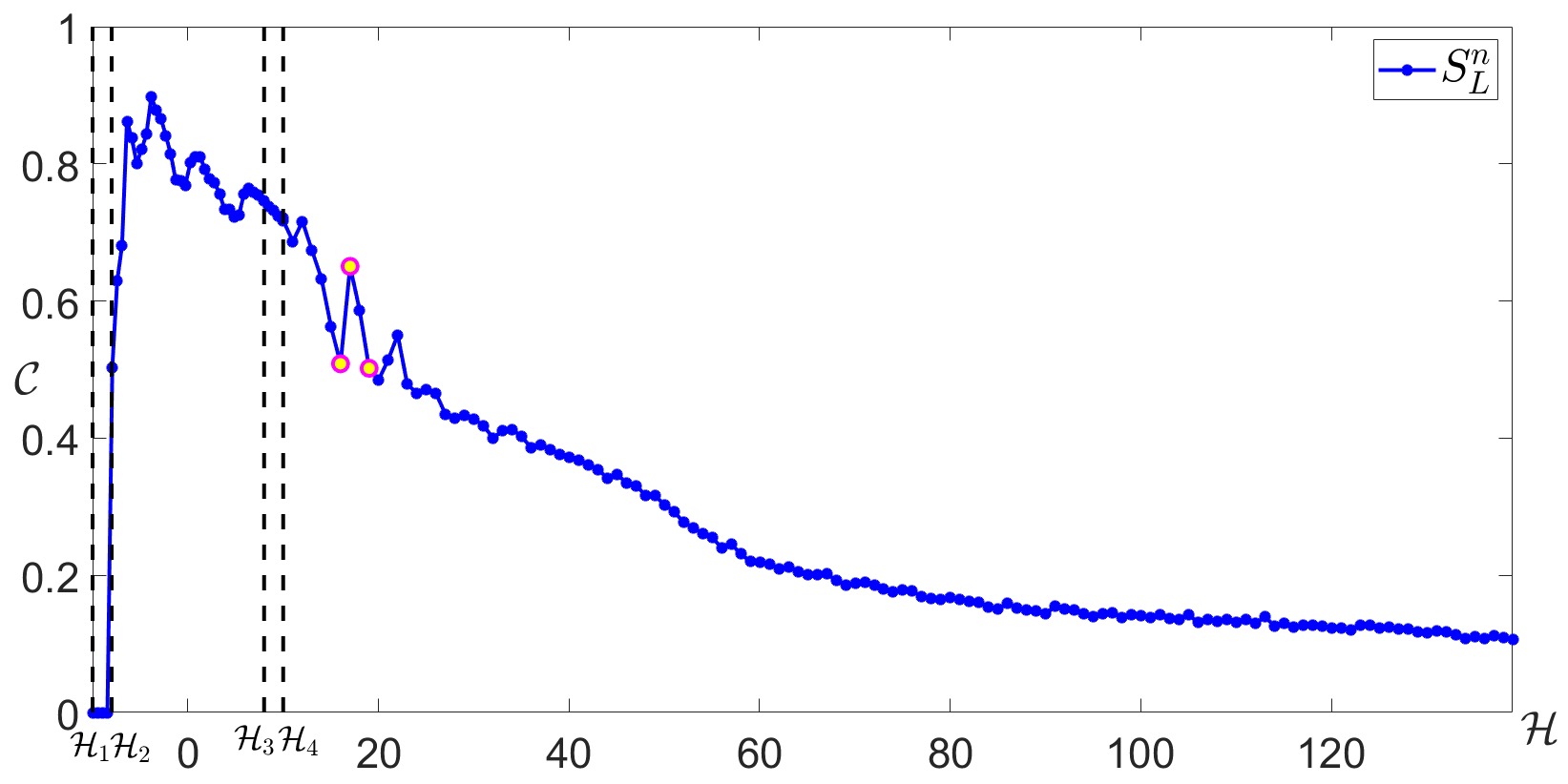} \\[.15cm]
    \includegraphics[scale=0.18]{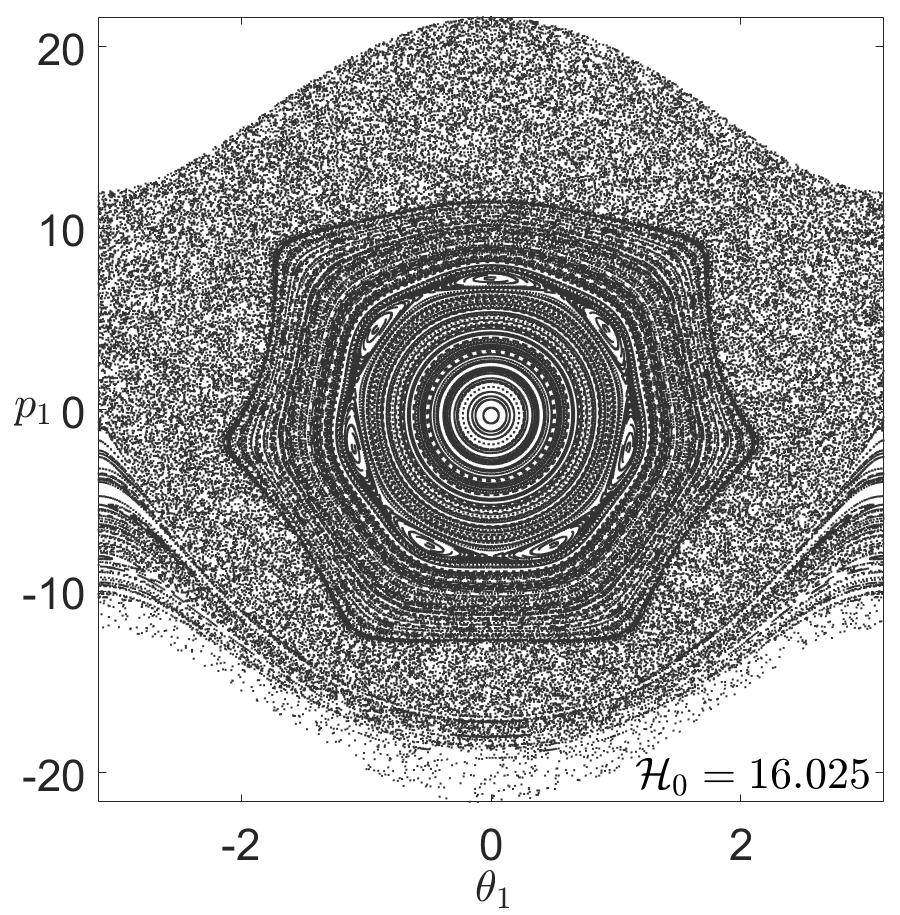}
    \includegraphics[scale=0.18]{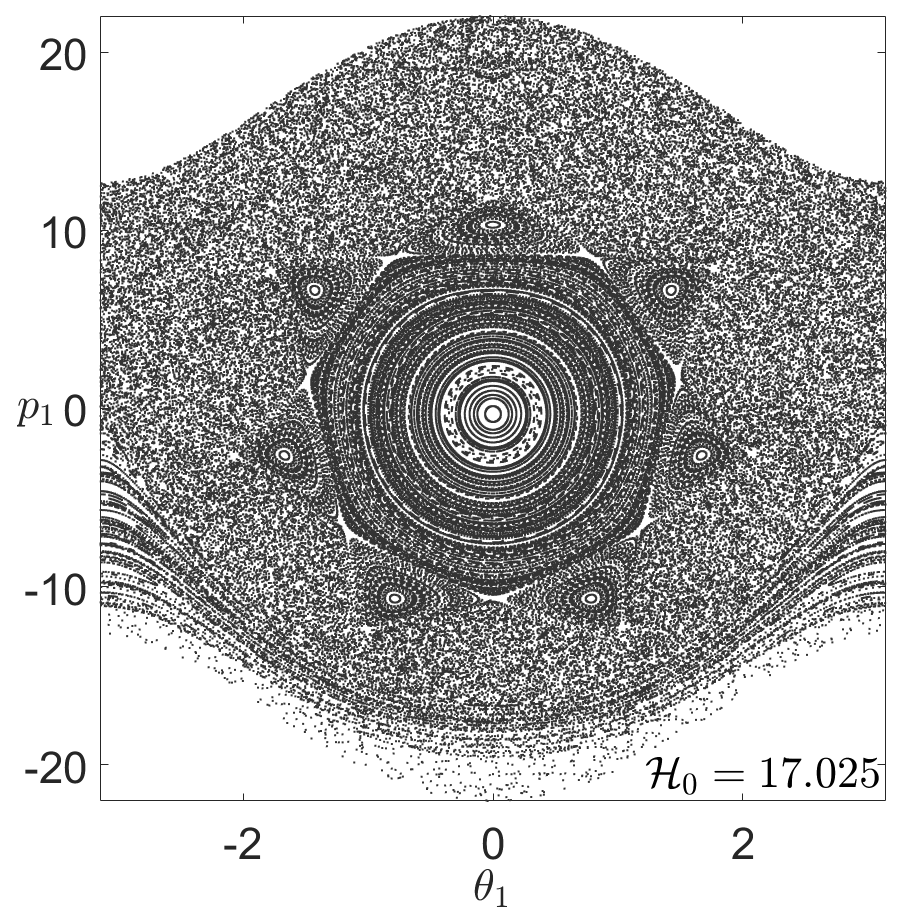}
    \includegraphics[scale=0.18]{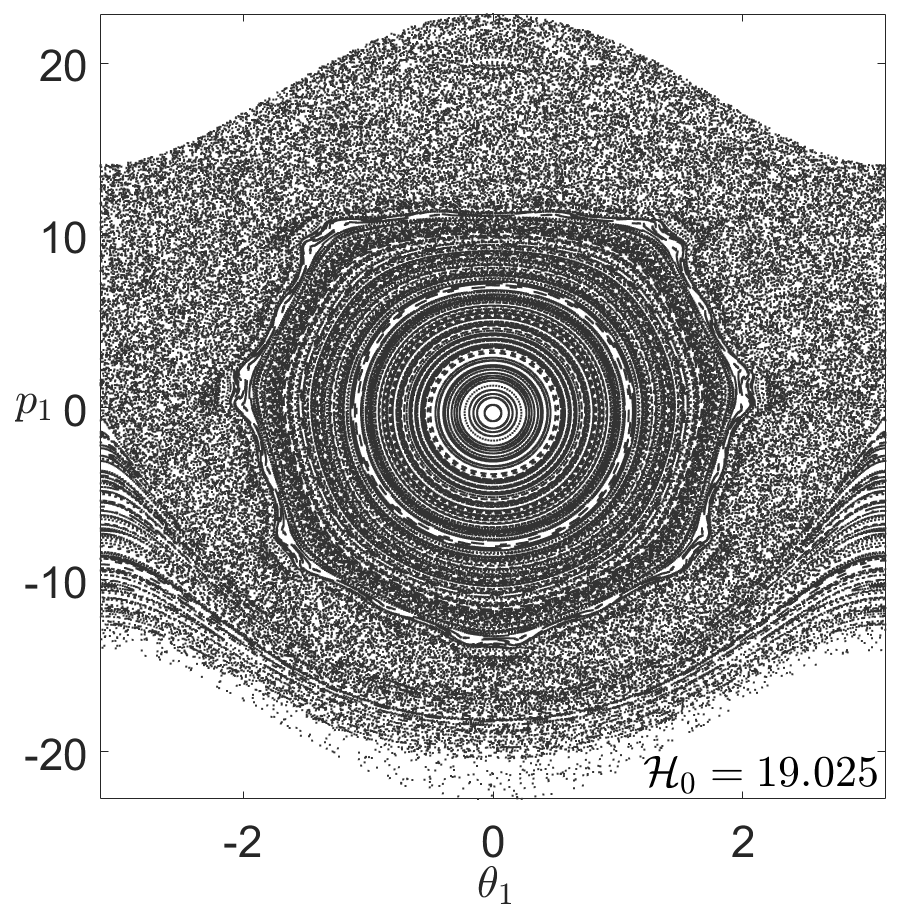}
    \caption{Chaotic fraction ($\mathcal{C}$) versus energy ($\mathcal{H}$) calculated with the chaos indicator $S^n_L$ for the model parameters $\alpha=1$ and $\sigma=8$. The black dotted lines are the energies of the local minimum ($\mathcal{H}_1$), of the index-1 saddles ($\mathcal{H}_2$ and $\mathcal{H}_3$) and that of the local maximum ($\mathcal{H}_4$) of the potential energy surface of the double pendulum. The panels on the bottom represent Poincaré sections for the energies marked on the blue curve above with yellow dots. A sudden peak in the chaotic fraction indicates a relevant change in the underlying phase space structure of the system.}
    \label{abrupt_chaos}
\end{figure}

Another interesting feature of the simulations we have carried out is that sudden peaks and valleys appear in the chaotic fraction of phase space trajectories. For example, in the case where $\alpha = 1$ and $\sigma = 8$, this behavior occurs in a short energy interval as displayed in Fig. \ref{abrupt_chaos}. As the system goes from an energy level of $\mathcal{H}_0 = 16.025$ to $\mathcal{H}_0 = 19.025$, the chaotic fraction increases and decreases abruptly by about $15\%$. This is a clear indication that significant changes are taking place in the phase space structure of the system due to bifurcations and the formation and destruction of KAM tori. These features are nicely captured by the Poincaré sections displayed below the chaotic fraction curve in Fig. \ref{abrupt_chaos}. Notice also that this behavior manifests repeatedly along the curve, but the peaks become less abrupt each time they occur.  

\begin{figure}[htbp]
    \centering
    A)\includegraphics[scale = 0.28]{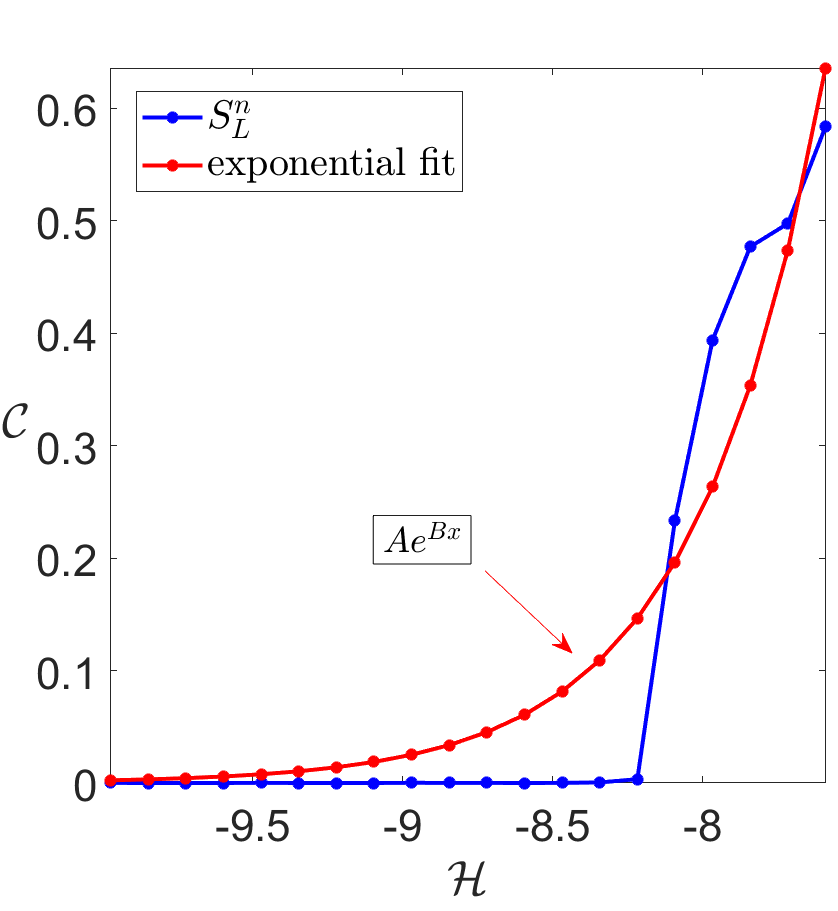}
    B)\includegraphics[scale = 0.28]{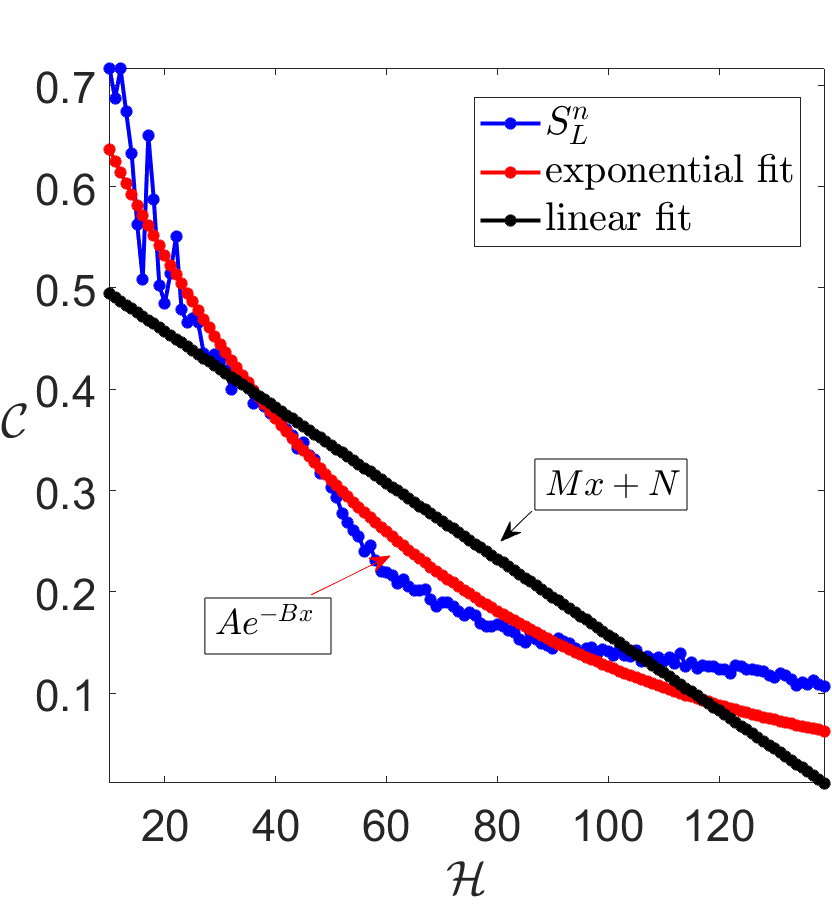} \\
    C) \hspace{.1cm} \includegraphics[scale = 0.9]{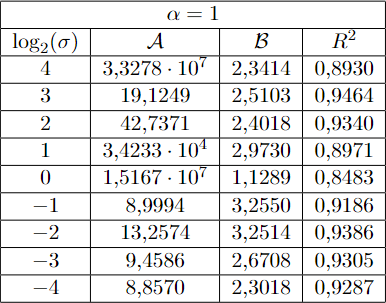} 
    D) \hspace{.1cm} \includegraphics[scale = 0.9]{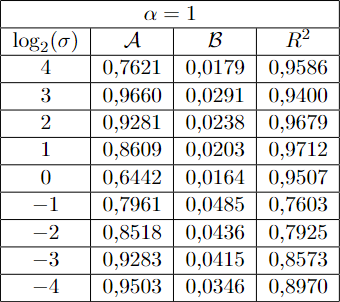}
    \caption{Fits performed for the chaotic fraction as a function of energy using exponential functions (red) and linear functions (black) for the model parameters $\alpha = 1$ and $\sigma = 8$. A) Exponential fit for the section of the curve in which the chaotic fraction grows; B) Exponential fit for the section of the graph where the chaotic fraction decreases; C) Coefficients of the exponential fit and $R^2$ indicator determined for different values of the model parameters in the section of the curve where the chaotic fraction increases; D) Coefficients of the exponential fit and $R^2$ indicator determined for different values of the model parameters in the energy regime where the chaotic fraction decreases.}
    \label{Exponential_fits}
\end{figure}

In addition to the study of the peaks that appear in the chaotic fraction curve, we have also performed different fits to study the trend of chaos growth and decay in the system with energy. Our simulations indicate that for values of the energy larger than that corresponding to the maximum of the potential, the curve behaves like a decreasing exponential of the form $\mathcal{A}e^{-\mathcal{B}x}$. For comparison, we have also considered a linear fit of the form $M x + N$, but the results obtained from this fit are much worse than those provided by the exponential model. This seems to indicate that chaos in the double pendulum system tends to follow an exponential decrease from the energy of the maximum. For example, if we take the case where $\alpha = 1$ and $\sigma = 8$, the $R^2$ coefficient for the exponential fit in this region is approximately $0.9586$, while for the linear fit the value obtained is $0.8048$, which shows that the linear fit is not suitable for this section of the curve as it is clearly shown in Fig. \ref{Exponential_fits} B). To support our claim, we have carried out this exponential fit for a fixed length ratio $\alpha = 1$ and different values of the mass ratio $\sigma$. A summary of the results obtained is included in the table displayed in Fig. \ref{Exponential_fits} D). It is important to note here that by using a fit with a decreasing exponential we have assumed that chaos tends to 0 in the system when the energy becomes arbitrarily large. However, we do not know if this assumption is true or not, even though the curve that represents the chaotic fraction as a function of the energy seems to point in that direction. Needless to say that one needs to be careful with these claims about the asymptotic behavior of the system. It is evident that one cannot make generalizations of this behavior to other dynamical systems. Indeed, there are many Hamiltonian systems in the literature where the chaotic fraction of the phase space does not decrease with energy, see e.g. \cite{Dahlqvist90}.

Another aspect that we have analyzed about the chaotic fraction as a function of the model parameters is the growth of chaos in the energy regime that covers the values from that of the minimum of the potential to the energy of the first index-1 saddle. For this section of the curve, we have also studied an exponential fit, since in this region the chaotic fraction grows very rapidly with energy. We have observed that an exponential function is perfectly suitable for this purpose, as can be seen in Fig. \ref{Exponential_fits} where the fit has the functional form $\mathcal{A}e^{\mathcal{B}x}$. We have performed a detailed study for the exponential models in the case of $\alpha = 1$ so that the values of the coefficients and the $R^2$ values for the fits are presented in the table included in Fig. \ref{Exponential_fits} C). Moreover, Fig. \ref{Exponential_fits} A) displays the section of the chaotic fraction curve where this growth occurs, and we have overlaid the exponential fit of the data in order to show how well the exponential fit performs in this approximation. These results are in the line of those discussed in \cite{cabrera2023regular}, where they address this behavior but only for the case where the pendulums have equal masses and lengths (that corresponds to our case $\alpha = 1$ and $\sigma = 0$).

\section{Conclusions} \label{Conclusions}

In this paper we have quantified the extent of chaos and regularity in the double pendulum by means of applying chaos indicators derived from LDs. Our findings show how these diagnostic tools can be used to effectively measure the degree of chaos in a computationally efficient manner. In order to study how the chaotic fraction of trajectories varies with energy and the other model parameters, we have shown that it is important to implement an adaptive strategy to set the threshold value for the chaos indicator that is required to accurately classify between chaotic and regular motion. The methodology is based on determining the location of the minimum between the two peaks that arise in the histogram of the chaos indicator obtained from the ensemble of initial conditions used in the simulations. This process is relevant since there is no universal threshold, as this critical value depends on many factors such as the integration time, the distance between neighbors, the model parameters and also the energy of the system.

Our parametric study of the degree of chaos in the double pendulum system has allowed us to establish how the masses and lengths of the pendulums influence the dynamics of this classical Hamiltonian system. The numerical analysis we have carried out reveals that 
the systems attains the maximum fraction of chaos when the pendulums have equal lengths. 
Another interesting feature of the problem that our simulations highlight is that, if the second pendulum is longer than the first one, chaos grows as the mass of the first pendulum is increased with respect to that of the second one. Conversely, if the first pendulum is longer, then chaos grows as the mass of the second pendulum is larger than the first one. 

Regarding the growth and decay of chaos in the system as a function of energy, we have 
found that in the low-energy regime, that is, from the energy at the minimum of the 
potential to that of the first index-1 saddle, chaos increases following an exponential law. This result was also reported in \cite{cabrera2023regular}, but only for one of the cases we cover in this work. By performing a curve fitting analysis, we have also observed from our numerical experiments that in the energy regime starting at approximately the value of the local maximum of the potential, the chaotic fraction shows an exponentially decreasing trend. This behavior seems to indicate that the system will recover full regularity in the high-energy regime. Interestingly, the reduction of chaos displayed by the system is not monotone, as abrupt peaks and valleys appear in the chaotic fraction curve. These sudden variations in the extent of chaos point in the direction that bifurcations and other significant changes are taking place in the underlying phase space of the system, as KAM islands are being generated and destroyed.

The work we have carried out in this paper opens new avenues for future research and raises many questions. For example, it would be interesting to analyze how chaos grows and decays in different Hamiltonian systems and, in particular, determine under what conditions full regularity is restored for the high-energy regime. Regarding the appearance of abrupt variations in the chaotic fraction, we believe that this phenomenon could be correlated with the complexity of the transitions that are taking in the phase space structure. It might be worth examining this relation by means of measuring the jaggedness of the chaotic fraction curve using for instance the fractal dimension. Another aspect of the problem that could be pursued is the development of control strategies on the system so that chaos varies smoothly with energy, reducing the sudden peaks and valleys that we have observed in the behavior of the double pendulum. Our research efforts in the future will focus on addressing all these questions, with the goal of gaining more insights and shedding some light on the intertwined chaotic and regular nature of Hamiltonian systems. We hope that this work could serve as a motivation for the further development of these chaos diagnostic techniques, as their simple implementation and computational efficiency make them ideal candidates to explore the chaotic behavior of any dynamical system, especially those that are high-dimensional and where the visualization of phase space structure is impractical.


\section*{Acknowledgements} \label{Acknowledgements}

The authors would like to acknowledge the support received for this research through the grant ``Introducci\'on a la Investigaci\'on'' from the Universidad de Alcal\'a.

%

\section*{Data availability}

The data that support the findings of this study are available from the corresponding author upon reasonable request.

\bibliography{referencias}

\end{document}